\documentclass[12pt]{article}%
\usepackage{graphicx}
\usepackage[intlimits]{amsmath}
\usepackage{latexsym}
\usepackage{amsfonts}
\usepackage{amssymb}%
 \makeatletter
\newfont{\footsc}{cmcsc10 at 8truept}
\newfont{\footbf}{cmbx10 at 8truept}
\newfont{\footrm}{cmr10 at 10truept}
\newtheorem{theorem}{Theorem}

\newtheorem{fact}[theorem]{Fact}

\newenvironment{proof}[1][Proof]{\noindent{\textbf {#1}  }}  {\hfill$\Box$\bigskip}

\begin{document}

\title{A spectral Erd\H{o}s-Stone-Bollob\'{a}s theorem}
\author{Vladimir Nikiforov\\{\small Department of Mathematical Sciences, University of Memphis, Memphis,
TN 38152}\\{\small email: vnikifrv@memphis.edu}}
\maketitle

\begin{abstract}
Let $r\geq3$ and $\left(  c/r^{r}\right)  ^{r}\ln n\geq1.$ If $G$ is a graph
of order $n$ and its largest eigenvalue $\mu\left(  G\right)  $ satisfies
\[
\mu\left(  G\right)  \geq\left(  1-\frac{1}{r-1}+c\right)  n,
\]
then $G$ contains a complete $r$-partite subgraph with $r-1$ parts of size
$\left\lfloor \left(  c/r^{r}\right)  ^{r}\ln n\right\rfloor $ and one part of
size greater than $n^{1-c^{r-1}}.$

This result implies the Erd\H{o}s-Stone-Bollob\'{a}s theorem, the essential
quantitative form of the Erd\H{o}s-Stone theorem.

Moreover, if $F$ is a fixed graph with chromatic number $r,$ then%
\[
\lim_{n\rightarrow\infty}\frac{1}{n}\max\left\{  \mu\left(  G\right)  :G\text{
is of order }n\text{ and }F\nsubseteq G\right\}  =1-\frac{1}{r-1}.
\]
This result implies the Erd\H{o}s-Stone-Simonovits theorem.\medskip

\textbf{Keywords: }\textit{largest eigenvalue; }$r$\textit{-partite subgraph;
Erd\H{o}s-Stone-Bollob\'{a}s theorem; Erd\H{o}s-Stone-Simonovits theorem}

\end{abstract}

This note is part of an ongoing project aiming to build extremal graph theory
on spectral grounds, see, e.g., \cite{BoNi07}, \cite{Nik02,Nik07h}.

The fundamental Erd\H{o}s-Stone theorem \cite{ErSt46} states that, given
$r\geq3$ and $c>0,$ every graph with $n$ vertices and $\left\lceil \left(
1-1/\left(  r-1\right)  +c\right)  n^{2}/2\right\rceil $ edges contains a
complete $r$-partite graph with each part of size $g\left(  n,r,c\right)  ,$
where $g\left(  n,r,c\right)  $ tends to infinity with $n$. In \cite{BES76}
Bollob\'{a}s and Erd\H{o}s found that $g\left(  r,c,n\right)  =\Theta\left(
\log n\right)  ,$ and in \cite{BoEr73}, \cite{BoKo94}, \cite{ChSz81}, and
\cite{Ish02} the function $g\left(  r,c,n\right)  $ was determined with great precision.

Here we deduce the Erd\H{o}s-Stone-Bollob\'{a}s result from a weaker, spectral condition.

Our notation follows \cite{Bol98}. Let $K_{r}\left(  s_{1},\ldots
,s_{r}\right)  $ be the complete $r$-partite graph with parts of size
$s_{1},\ldots,s_{r},$ and let $\mu\left(  G\right)  $ be the largest adjacency
eigenvalue of a graph $G$. Our main result:

\begin{theorem}
\label{genZ}Let $r\geq3,$ $\left(  c/r^{r}\right)  ^{r}\ln n\geq1,$ and $G$ be
a graph with $n$ vertices. If
\[
\mu\left(  G\right)  \geq\left(  1-\frac{1}{r-1}+c\right)  n,
\]
then $G$ contains a $K_{r}\left(  s,\ldots s,t\right)  $ with $s\geq
\left\lfloor \left(  c/r^{r}\right)  ^{r}\ln n\right\rfloor $
and$\ t>n^{1-c^{r-1}}.$
\end{theorem}

As an easy consequence, we strengthen the Erd\H{o}s-Stone-Simonovits theorem
\cite{ErSi66}.

\begin{theorem}
\label{ErSi}Let $r\geq3$ and $F$ be a fixed graph with chromatic number $r.$
Then%
\[
\lim_{n\rightarrow\infty}\frac{1}{n}\max\left\{  \mu\left(  G\right)  :G\text{
is of order }n\text{ and }F\nsubseteq G\right\}  =1-\frac{1}{r-1}.
\]

\end{theorem}

\subsubsection*{Remarks}

\begin{itemize}
\item[-] Since $\mu\left(  G\right)  $ is at least the average degree of $G,$
Theorem \ref{genZ} implies the following form of the
Erd\H{o}s-Stone-Bollob\'{a}s theorem: \medskip

\emph{Let }$r\geq2,$\emph{ }$\left(  c/r^{r}\right)  ^{r}\ln n\geq1,$\emph{
and }$G$\emph{ be a graph with }$n$\emph{ vertices. If }$G$\emph{ has
}$\left\lceil \left(  1-1/\left(  r-1\right)  +c\right)  n^{2}/2\right\rceil
$\emph{ edges, then }$G$\emph{ contains a }$K_{r}\left(  s,\ldots,s,t\right)
$\emph{ with }$s\geq\left\lfloor \left(  c/r^{r}\right)  ^{r}\ln
n\right\rfloor $\emph{ and\ }$t>n^{1-c^{r-1}}.$\medskip

This is slightly stronger than the result in \cite{BoEr73} and is comparable
to the results in \cite{BoKo94}.

\item[-] The relation between $c$ and $n$ in Theorem \ref{genZ} needs
explanation. First, for fixed $c,$ it shows how large must be $n$ to get a
valid conclusion. But, in fact, the relation is subtler, for $c$ itself may
depend on $n,$ e.g., letting $c=1/\ln\ln n,$ the conclusion is meaningful for
sufficiently large $n.$

\item[-] Using random graphs, we see that almost all graphs on $n$ vertices
contain no $K_{2}\left(  s,s\right)  $ with $s$ larger than $C\log n$ for some
$C>0,$ independent of $n.$ Hence, for constant $c,$ Theorem \ref{genZ} is
essentially best possible.

\item[-] After this note has been made public, we were informed that Babai and
Guiduli \cite{Gui96} have proved Theorem \ref{ErSi} using the Szemer\'{e}di
Regularity Lemma; for a recent account on this matter see \cite{BaGu07}.
\end{itemize}

\subsubsection*{Proofs}

A word about our proof methods seems in place. In an ongoing series of papers,
e.g., \cite{BoNi07,BoNi07a}, \cite{Nik02,Nik07h}, we are developing a set of
wide-range tools for use in extremal and spectral graph theory. Sometimes, as
seen below, these tools are so effective that the proofs become vanishingly short.

Write $k_{r}\left(  G\right)  $ for the number of $r$-cliques of a graph $G$.
The following facts play crucial roles in our proof of Theorem \ref{genZ}.

\begin{fact}
[\cite{BoNi07}, Theorem 2]\label{leNSMM}If $r\geq2$ and $G$ is a graph of
order $n,$ then%
\[
k_{r}\left(  G\right)  \geq\left(  \frac{\mu\left(  G\right)  }{n}-1+\frac
{1}{r}\right)  \frac{r\left(  r-1\right)  }{r+1}\left(  \frac{n}{r}\right)
^{r}.
\]
$\hfill\square$
\end{fact}

\begin{fact}
[\cite{Nik07a}, Theorem 1]\label{ES}Let $r\geq2,$ $c^{r}\ln n\geq1,$ and $G$
be a graph of order $n$. If $k_{r}\left(  G\right)  \geq cn^{r},$ then $G$
contains a $K_{r}\left(  s,\ldots s,t\right)  $ with $s=\left\lfloor c^{r}\ln
n\right\rfloor $ and $t>n^{1-c^{r-1}}.\hfill\square$
\end{fact}

\begin{fact}
\label{tsize}The number of edges of $T_{r}\left(  n\right)  $ satisfies
$2e\left(  T_{r}\left(  n\right)  \right)  \geq\left(  1-1/r\right)
n^{2}-r/4.\hfill\square$
\end{fact}

\begin{proof}
[\textbf{Proof of Theorem \ref{genZ}}]In view of $\mu\left(  G\right)
\geq\left(  1-1/\left(  r-1\right)  +c\right)  n,$ Fact \ref{leNSMM} implies
that
\[
k_{r}\left(  G\right)  >c\frac{r-2}{r^{r}}n^{r}\geq\frac{c}{r^{r}}n^{r}.
\]
Hence, Fact \ref{ES} implies that $G$ contains a $K_{r}\left(  s,\ldots
,s,t\right)  $ with%
\[
s\geq\left(  c/r^{r}\right)  ^{r}\ln n,\text{ \ }t>n^{1-c^{r-1}},
\]
completing the proof.
\end{proof}

\bigskip

\begin{proof}
[\textbf{Proof of Theorem \ref{ErSi}}]Theorem \ref{genZ} implies that
\[
\limsup_{n\rightarrow\infty}\frac{1}{n}\max\left\{  \mu\left(  G\right)
:G\text{ is of order }n\text{ and }F\nsubseteq G\right\}  \leq1-\frac{1}%
{r-1}.
\]
On the other hand, writing $T_{s}\left(  n\right)  $ for the $s$-partite
Tur\'{a}n graph of order $n,$ in view of Fact \ref{tsize}, we see that
\[
\frac{\mu\left(  T_{r-1}\left(  n\right)  \right)  }{n}\geq\frac{2e\left(
T_{r-1}\left(  n\right)  \right)  }{n^{2}}\geq1-\frac{1}{r-1}-\frac
{r-1}{4n^{2}}.
\]
Since $T_{r-1}\left(  n\right)  $ is $\left(  r-1\right)  $-partite, it
contains no copy of $F$. Therefore,
\[
\liminf_{n\rightarrow\infty}\frac{1}{n}\max\left\{  \mu\left(  G\right)
:G\text{ is of order }n\text{ and }F\nsubseteq G\right\}  \geq1-\frac{1}%
{r-1},
\]
completing the proof.
\end{proof}

\subsubsection*{Concluding remark}

Finally, a note about the project mentioned in the introduction: in this
project we aim to give wide-range results that can be used further, adding
more integrity to spectral extremal graph theory.\bigskip

\textbf{Acknowledgement} Thanks Felix Lazebnik for the motivation to write
this note, and L\'{a}szl\'{o} Babai for useful suggestions.

\end{document}